\documentclass[12pt,a4paper]{article}

\usepackage{amsfonts,latexsym}

\makeatletter

\@addtoreset{equation}{section}
\newcommand{\l@abcd}[2]{\hbox to\textwidth{#1\dotfill #2}}

\makeatother

\newcommand*{\mybegintheorem}[1]{\begin{trivlist}\it%
      \item[\hspace{\labelsep}{\bf #1}]}
\newcommand*{\myendtheorem}{\end{trivlist}}
\newenvironment*{theorem*}{\mybegintheorem{Theorem.}}{\myendtheorem}
\newenvironment*{proposition*}{\mybegintheorem{Proposition.}}{\myendtheorem}
\newenvironment*{lemma*}{\mybegintheorem{Lemma.}}{\myendtheorem}
\newenvironment*{corollary*}{\mybegintheorem{Corollary.}}{\myendtheorem}
\newenvironment*{definition*}{\mybegintheorem{Definition.}}{\myendtheorem}

\renewcommand{\phi}{\varphi}
\renewcommand{\epsilon}{\varepsilon}

\newcommand{\lra}{\longrightarrow}

\newcommand{\PTw}{{{\mathbb P}^2}}

\newcommand{\PP}{{{\mathbb P}}}
\newcommand{\ZA}{{\mathbb Z}}
\newcommand{\QA}{{\mathbb Q}}
\newcommand{\RA}{{\mathbb R}}
\newcommand{\CA}{{\mathbb C}}

\title{Normal surfaces with strictly nef anticanonical divisors}
\author{Mikhail Grinenko}
\date{{\small Steklov Institute, Moscow\\e-mail: grin@mi.ras.ru}}

\begin{document}
\maketitle

\paragraph{Abstract.} {\it\small In this paper we study normal surfaces
whose anticanonical divisors are strictly nef, i.e. $(-K)\cdot C>0$
for every curve $C$.}

\paragraph{\bf Remark.} When this paper was ready to print
Keiji Oguiso pointed me on the paper of Fernando Serrano [6]
where some results of my work already proved (I mean two-dimensional case).
Nonetheless in this paper we work in more general case.

\paragraph{}
It is well-known that a lot of properties of algebraic
varieties are defined by the numerical properties of their
anticanonical divisors. We consider algebraic surfaces from this point
of view.

\paragraph{1. Smooth surfaces.} Let $X$ be a smooth projective surfaces
over the field of complex numbers $\CA$. Then $X$ is rational if and only if
$$
    H^0(X,{\cal O}(2K_X)) = H^1(X,{\cal O}(K_X)) = 0.
$$
These conditions are always met when $-K_X$ is ample, by the Kodaira
vanishing theorem. On the other hand, the Nakai-Moishezon criterion yields
$$
 \mbox{ Divisor $H$ is ample } \Longleftrightarrow
 \mbox{ $H^2>0$ and $H\cdot C>0$ for every curve C}
$$
As it has been shown in [3], we may drop the condition $H^2>0$
when $H$ is an anticanonical divisor. One can get another description
of this result.
\begin{itemize}
\item[]{\it Let ${\mathbb {NE}}(X)$ be a cone of effective 1-cycles on X.
Then $X$ is Del Pezzo if and only if $-K_X$ is positive on
${\mathbb {NE}}(X)$.}
\end{itemize}

\paragraph{2. Normal surfaces.} One may try to generalize the previous result.
Let $X$ be a normal surface. As in [5], we can take a canonical Weil
divisor which define an invertible sheaf on the nonsingular locus
of $X$. Because the intersection theory on normal surfaces does exist,
let $-K_X\cdot C>0$ for every curve $C\subset X$. What is possible
to say in this situation? Let's consider some examples.

\paragraph{Examples.} {\bf a)} Let $Y$ be a rational ruled surface with
$C_0$ as an exceptional section, $C_0^2\le -2$. If $\phi: Y\to X$ is
a contraction of $C_0$ then $X$ is rational and $\QA$-Gorenstein
([1], [5]). Moreover, $-mK_X$ is ample for some integer
$m > 0$.

\noindent{\bf b)} Let $C$ be a smooth curve of genus 2 and
${\cal L}={\cal O}_C(E)$, where $E$ is a divisor of degree
$\deg E=-3$ on $C$. Let's consider a ruled surface
$$
     \pi: Y=\PP_C({\cal E})\to C
$$
for the sheaf ${\cal E}={\cal O}_C\oplus E$ of rank 2. Thus
$$
   K_Y\sim -2C_0+\pi^*(F+E)
$$
where $C_0$ is an exceptional section, $C_0^2=-3$, and F is a canonical
divisor on C, $\deg F=2$.

If $\phi: Y\to X$ is a contraction of $C_0$ to the normal surface $X$,
we have
$$
   3\phi^*K_X\sim -C_0+3\pi^*(F+E),
$$
so $-K_X$ is numerically ample on X (i.e. $(-K_X)^2>0$ and $(-K_X)\cdot C>0$
for every curve $C$). It is easy to see that $X$ is projective. Indeed,
let's take a divisor $H\sim C_0-\pi^*E$ on $Y$. Then $H\cdot C_0=0$,
$H\cdot C>0$ for every curve $C\ne C_0$, $H^2 >0$ and
${\cal O}_Y(H)\otimes{\cal O}_{C_0}\simeq{\cal O}_{C_0}$, so $H=\phi^*D$
for some ample Cartier divisor $H$ on $X$ ([5]).

Since $(-C_0+3\pi^*(F+E))|_{C_0}\sim 3F+2E$, $X$ is
$\QA$-Gorenstein if and only if $n(3F+2E)\sim 0$ for some integer $n\ne 0$.
Anyway, $X$ is not rational.

\paragraph{3.} Now I formulate the main result of this paper.

\begin{theorem*}
Let $X$ be a normal projective surface over the field $\CA$ such that
$(-K_X)\cdot C>0$ for every curve $C$ on $X$. Then
\begin{itemize}
\item[(i)] $(-K_X)^2 > 0$;
\item[(ii)] $X$ is rational if and only if all singularities of $X$ are
rational, i.e. $R^1\pi_*{\cal O}_Y=0$ for some resolution $\pi:Y\to X$;
\item[(iii)] $X$ is rational if and only if $X$ is $\QA$-factorial.
\end{itemize}
\end{theorem*}

\paragraph{4.} {\bf Proof.} {\it(i)} We need the following
lemma:

\begin{lemma*}
Let $X$ be a surface with Du Val singulariries, and $(-K_X)\cdot C>0$
for every curve $C$. Then $(-K_X)^2>0$.
\end{lemma*}
\noindent{\it Proof of lemma.} By [3], we can suppose $X$ to be singular.
If $(-K_X)^2=0$, then $|nK_X|=\emptyset$ for any $n>0$. Let
$$
   f: Y\lra X
$$
be a minimal resolution of $X$. We have $K_Y=f^*K_X$ and
$k(X)=k(Y)=-\infty$. Moreover, from the Riemann-Roch theorem
$$
    -h^1(Y,-K_Y)=1-q(Y) ,
$$
where $q(Y)$ is the irregularity of $Y$. Then $q(Y)\ge 1$, and the minimal
model of $Y$ is a ruled surface over the curve of genus $g\ge 1$. But
$$
   0=K_Y^2\le 8(1-g)\le 0 ,
$$
so $g=1$ and $Y$ is minimal. It is impossible because $Y$ contains the
smooth rational curve $E$ with $E^2=-2$. This contradiction proves the lemma.

In order to prove the theorem, let's take the minimal resolution
$\pi: Y\to X$, so
$$
    K_Y=\pi^*K_X+\sum\alpha_iE_i ,
$$
where $\{E_i\}$ are exceptional curves and all $\alpha_i\le 0$. By the
previous lemma, we can suppose $\alpha_i\ne 0$ for some $i$.

It is easy to see that $K_Y\cdot E_i\ge 0$ for all $i$ and
$$
    K_Y\cdot C\le K_X\cdot\pi_*C < 0
$$
for any curve $C\notin\{E_i\}$.

By the Cone theorem,
$$
   \overline{{\mathbb{NE}}}(Y)=%
   \overline{\left\{\overline{\mathbb{NE}}_{>0}+\sum R_j\right\}} ,
$$
where $\{R_j\}$ are extremal rays on $Y$. Let
$K_Y^{\perp}=\{z\in\overline{{\mathbb{NE}}}(Y): z\cdot K_Y=0\}$.
Then
$$
   \overline{{\mathbb{NE}}}(Y)=%
   \overline{\left\{\sum\RA_{+}[E_i]+K_Y^{\perp}+\sum R_j\right\}} .
$$
Let the extremal ray $R_j$ be generated by the class of the curve $C_j$,
i.e. $R_j=\RA_{+}[C_j]$. Since $Y\ne\PTw$, we have either $C_j$ is an
exceptional curve of the first kind, or $Y$ is a ruled surface with $C_j$ as
its fiber ([4]). In the lather case $K_Y=\pi^*K_X-qE$, $q>0$, $E$ is
an irreducible curve with $E^2=-e < 0$. Let $F$ be a fiber of $Y$, and
$(aE+bF)$ be a $\QA$-divisor such that $(aE+bF)\sim -\pi^*K_X$. We have
$$
\begin{array}{l}
(aE+bF)\cdot F=a > 0 \\
(aE+bF)\cdot E= -ae+b=0
\end{array}
$$
hence $K_X^2=(aE+bF)^2=a^2e > 0$, and the theorem is proved in this case.

So we will suppose $C_j$ to be a {\bf-1}-curve, $K_Y\cdot C_j=-1$. Let
$E=\sum E_i$ and  $mE$ be a Cartier divisor for some integer $m>0$. Then
$D=-mK_Y+mE$ is Cartier, $D\equiv m\pi^*(-K_X)$. Since
$$
   m=D\cdot C_j-mE\cdot C_j ,
$$
$-mE$ is effective and $C_j\notin\{E_i\}$, we have
$$
\begin{array}{l}
  \min_{j}\{D\cdot C_j\}\ge 1 \\
  \max_{j}\{(-mE)\cdot C_j\}\le m-1
\end{array}
$$
Let's consider the function $\bar D=D+E$ on the cone
$\overline{{\mathbb{NE}}}(Y)$. Then
\begin{itemize}
\item[(1)] $\bar D\ge 0$ on the cone $\sum\RA_{+}[E_i]$ because
   $\bar D\cdot E_i=E\cdot E_i\ge 0$;
\item[(2)] $\bar D\ge 0$ on the cone $K_Y^{\perp}$: if $z\in K_Y^{\perp}$
     then $z\cdot E=-\pi^*K_X\cdot z\ge 0$;
\item[(3)] $\bar D\ge 0$ on the cone $\sum R_j$ since
$D\cdot C_j=D\cdot C_j+E\cdot C_j\ge 1-\frac{m-1}{m}>0$.
\end{itemize}
Thus $\bar D\ge 0$ on $\overline{{\mathbb{NE}}}(Y)$, hence
$\bar D^2=m^2K_X^2+E^2\ge 0$, and then $K_X^2 > 0$ because $E^2 < 0$.
The part {\it(i)} is proved.

{\it(ii)} Since $(-K_X)$ is {\it nef} and {\it big},
we have $H^i(X,{\cal O}_X)=0$ for $i>0$ ([5], theorem 5.1).
Let $\pi:Y\to X$ be a resolution of singularities,
$\pi_*{\cal O}_Y={\cal O}_X$.

Let's suppose that all singularities are rational, $R^1\pi_*{\cal O}_Y=0$.
Then
$$
   H^i(Y,{\cal O}_Y)=H^i(X,\pi_*{\cal O}_Y)=H^i(X,{\cal O}_X)=0
$$
for all i. Moreover, $H^0(Y,2K_Y)=0$, thus both $Y$ and $X$ are rational.

Conversely,  $H^1(Y,{\cal O}_Y)=0$ if $X$ is rational.  We have the
spectral sequence
$$
   E_2^{p,q}=H^p(X,R^q\pi_*{cal O}_Y)\Rightarrow H^{p+q}(Y,{\cal O}_Y)
$$
which yields the next exact sequence:
$$
   0\lra H^1(X,{\cal O}_X)\lra H^1(Y,{\cal O}_Y)\lra %
    H^0(X,R^1\pi_*{\cal O}_Y)\lra H^2(X,{\cal O}_X),
$$
thus $R^1\pi{\cal O}_Y=0$.

\noindent{\it (iii)} Let's take a resolution of singularities
$$
   \pi: Y\lra X
$$
with exceptional curves $E_1,\ldots,E_s$, and
$$
  E=\sum_{i=1}^sE_i
$$

If $X$ is rational, then $R^1\pi_*{\cal O}_Y=0$, hence
$R^1\pi_*{\cal O}_E=0$, i.e. $H^1(E,{\cal O}_E)=0$. By [1],
$$
    H^1(E,{\cal O}^*_E)\simeq \ZA^s
$$
in this case. For any Weil divisor $D$ on $X$ we have
$$
   \pi^*{\cal O}_X(mD)\in Pic(Y)
$$
for some integer $m>0$, and $\deg_{E_i}{\cal O}_E(\pi^*mD)=0$ for every i,
so
$$
      {\cal O}_E(\pi^*mD)=(0,\ldots,0)\in\ZA^s.
$$
Then ${\cal O}_E(\pi^*mD)\otimes{\cal O}_E\simeq{\cal O}_E$ and $mD$ is
Cartier divisor on $X$ ([5]).

Conversely, if $X$ is not rational, then
$$
   \chi({\cal O}_X)=1\ne 1-h^1(Y,{\cal O}_Y)=\chi({\cal O}_Y),
$$
and there is a nonrational curve among $E_1,\ldots,E_s$
(for example $E_1$). Let
$$
    f: Y\lra Z
$$
be a birational morphism to the minimal model of $Y$, and $C=f(E_1)$.
$Z$ is ruled, and we can choose two fibers $F_1$ and $F_2$ such that
$f$ is an isomorphism near them and
$$
     {\cal O}_Z(n(F_1-F_2))\otimes{\cal O}_C\ne{\cal O}_C
$$
for every integer $n\ne 0$. Then $\pi(f^{-1}(F_1)-f^{-1}(F_2))$ is a Weil
divisor on $X$ which is not $\QA$-Cartier.

The proof is complete.

\paragraph{6. Remarks.}{\it 1)} The part {\it(ii)} of the theorem works
only under the condition on $(-K_X)$ to be {\it big} and {\it nef}. It is
easy to get a counter-example. We can blow up 12 points on $\PTw$ which
are an intersection of the smooth elliptic curve and some quartic. The
contraction of the strict transform of this elliptic curve yields
a normal projective rational surface with a nonrational singurarity.

\noindent{\it 2)} A normal surface with a numerically ample anticanonical
divisor is called {\it numerical Del Pezzo}. Such surfaces with
nonrational singularities are described in [2].

\noindent{\it 3)} We can use the same reasons to obtain the next result
in {\it log}-theory:

\begin{proposition*}
Let $(X,B)$ be a pair with log terminal singularities. The $\QA$-Cartier
divisor $-(K_X+B)$ is ample if and only if it is strictly nef, i.e.
$-(K_X+B)\cdot C >0$ for every curve $C$.
\end{proposition*}

\paragraph{\Large References}

\paragraph{[1]} Artin M. {\it Some numerical criteria for contractibility of
curves on algebraic surfaces,}  Amer. J. Math. {\bf 84} (1962), 485-496

\paragraph{[2]} Cheltsov I.  {\it Del Pezzo surfaces with nonrational
singularities,} Mathemati\-cheskie Zametki, {\bf 62} (1997),
451-467 (in russian)

\paragraph{[3]} Maeda H.  {\it A criterion for a smooth surface to be Del
Pezzo,} Math. Proc. Cambridge Phil. Soc. {\bf 113} (1993), 1-3

\paragraph{[4]} Mori S.  {\it Threefolds whose canonical bundle are not
numerically effective,} Ann. of Math. {\bf 116} (1982), 283-360

\paragraph{[5]} Sakai F.  {\it Weil divisors on normal surfaces,}
Duke Math. J. {\bf 51} (1984), 877-887

\paragraph{[6]} Serrano F.  {\it Strictly nef divisors and Fano threefolds,}
J. Reine Angew. Math. {\bf 464} (1995), 187-206

\end{document}